\documentclass[12pt,oneside,letterpaper]{article}%
\usepackage[centertags]{amsmath}
\usepackage[latin1]{inputenc}
\usepackage{amsfonts}
\usepackage{amssymb}
\usepackage{graphicx}
\usepackage{acronym}%
\newcommand{\be}{\beta}
\newcommand{\ds}{\displaystyle }
\newcommand{\f}{\frac}
\newcommand{\p}{\partial}
\newcommand{\bb}{\begin{equation}}
\newcommand{\e}{\end{equation}}
\newcommand{\ba}{\begin{array}}
\newcommand{\ea}{\end{array}}
\def\b<#1>{\langle #1 \rangle}
\newtheorem{theorem}{Theorem}

\newtheorem{corollary}[theorem]{Corollary}

\newtheorem{lemma}[theorem]{Lemma}

\newenvironment{proof}[1][Proof]{\noindent\textbf{#1.} }{\ \rule{0.5em}{0.5em}}
\begin{document}

\title{Noether Symmetries and Conservation{\huge \textbf{ }}Laws For Non-Critical
Kohn-Laplace{\huge \textbf{ }}Equations on\textbf{ }%
Three-Dimensional{\huge \textbf{ }}Heisenberg Group}
\author{\textrm{{\large Igor Leite Freire}}\\Institute of Mathematics, Statistics and Scientific Computation\\IMECC-UNICAMP CP 6065\\13083-859 Campinas, SP, Brazil\\e-mail: igor@ime.unicamp.br \\{\footnotesize 2000 AMS Mathematics}{ }{\footnotesize Classification
numbers:35H10, 58J70}}
\date{\ }
\maketitle

\begin{abstract}
We show which Lie point symmetries of non-critical semilinear Kohn-Laplace
equations on the Heisenberg group $H^{1}$ are Noether symmetries and we
establish their respectives conservations laws.

\end{abstract}

\newpage

\section{Introduction and Main Results}

\ 

In this paper we show which Lie point symmetries of the semilinear Kohn -
Laplace equations on the three-dimensional Heisenberg group $H^{1}$,
\begin{equation}
\label{klgen}\Delta_{H^{1}} u+f(u)=0,
\end{equation}
are Noether's symmetries, and we establish their respectives conservation laws.

The Kohn - Laplace operator on $H^{1}$ is defined by
\[
\Delta_{H^{1}}:=X^{2}+Y^{2}=\frac{\partial^{2}}{\partial x^{2}}+\frac
{\partial^{2}}{\partial y^{2}}+4(x^{2}+y^{2})\frac{\partial^{2}}{\partial
t^{2}}+4y\frac{\partial^{2}}{\partial x\partial t}-4x\frac{\partial^{2}%
}{\partial y\partial t},
\]
where%
\[
X=\frac{\partial}{\partial x}+2y\frac{\partial}{\partial t}\text{, }%
Y=\frac{\partial}{\partial y}+2x\frac{\partial}{\partial t}.
\]
\ \ 

Eq.(\ref{klgen}) possesses variational structure and can be derived from the
Lagragian
\begin{equation}
\mathcal{L}= \frac{1}{2}u_{x}^{2}+\frac{1}{2}u_{y}^{2}+2(x^{2}+y^{2})u_{t}
^{2}+2yu_{x}u_{t}-2xu_{y}u_{t}-F(u),\label{lag}%
\end{equation}
with $F^{\prime}(u)=f(u)$.

The group structure, the left invariant vector fields on $H^{1}$ and their Lie
algebra are given, respectively, by $\phi:\mathbb{R}^{3}\times\mathbb{R}%
^{3}\rightarrow\mathbb{R}^{3}$, where
\[
\phi((x,y,t),(x_{0},y_{0},t_{0})):=(x+x_{0},y+y_{0},t+t_{0}+2(xy_{0}-yx_{0})),
\]%
\begin{equation}\label{fields}
\ba{l c l}
X  & = &\ds{\frac{d}{ds}\phi((x,y,t),(s,0,0))|_{s=0}{=}\frac{\partial}{\partial
x}+2y\frac{\partial}{\partial t},} \\ \\
Y  & = &\ds{ \frac{d}{ds}\phi((x,y,t),(0,s,0))|_{s=0}{=}\frac{\partial}{\partial
y}+2x\frac{\partial}{\partial t}},\\ \\
Z  & = & \ds{\frac{d}{ds}\phi((x,y,t),(0,0,s))|_{s=0}{=}\frac{\partial}{\partial
t}},\ea
\end{equation}
and
\[
\lbrack X,T]=[Y,T]=0,\;\;[X,Y]=-4T.
\]

In \cite{gc} a complete group classification for equation (\ref{klgen}) is
presented. It can be summarized as follows.

Let $G_{f}:=\{T,R,\tilde{X},\tilde{Y}\}$, where
\begin{equation}
T=\frac{\partial}{\partial t},\;\;R=y\frac{\partial}{\partial x}%
-x\frac{\partial}{\partial y},\;\;\tilde{X}=\frac{\partial}{\partial
x}-2y\frac{\partial}{\partial t},\text{ and }\tilde{Y}=\frac{\partial
}{\partial y}+2x\frac{\partial}{\partial t}.\label{kil}%
\end{equation}
For any function $f(u)$, the group $G_{f}$ is a (sub)group of symmetries. Its
Lie algebra is summarized in Table \ref{tablegeneral}. \begin{table}[h]
\label{tablegeneral}
\par
\begin{center}%
\begin{tabular}
[c]{|c|c|c|c|c|}\hline
& T & R & $\tilde{X}$ & $\tilde{Y}$\\\hline
T & $0$ & $0$ & $0$ & $0$\\\hline
R & $0$ & $0$ & $\tilde{Y}$ & -$\tilde{X}$\\\hline
$\tilde{X}$ & $0$ & -$\tilde{Y}$ & $0$ & $4$T\\\hline
$\tilde{Y} $ & $0$ & $\tilde{X}$ & - $4$T & $0$\\\hline
\end{tabular}
\end{center}
\caption{{\small {Lie brackets of equation $(\ref{klgen})$ with $f(u)$}
arbitrary.}}%
\end{table}

For special choices of function $f(u)$ in (\ref{klgen}), the symmetry group
can be enlarged. Below we exhibit these functions and their respective
additional symmetries and Lie algebras.

\begin{itemize}
\item If $f(u)=0$, the additional symmetries are
\begin{equation}
\label{V1}\ba{l c l} V_{1} & = & \ds{(xt-x^{2}y-y^{3})\frac{\partial}{\partial x} +
(yt+x^{3}+xy^{2})\frac{\partial}{\partial y}} \\ \\
& & \ds{+ (t^{2}-(x^{2}+y^{2})^{2}%
)\frac{\partial}{\partial t}- t u \frac{\partial}{\partial u}},
\ea\end{equation}
\begin{equation}
\label{V2}\ba{ l c l} V_{2} & = & \ds{(t-4xy)\frac{\partial}{\partial x} + (3x^{2}-y^{2}%
)\frac{\partial}{\partial y}}\\ \\
& & \ds{ - (2yt+2x^{3}+2xy^{2})\frac{\partial}{\partial t}
+ 2 y u \frac{\partial}{\partial u}},
\ea\end{equation}
\begin{equation}
\label{V3}\ba{ l c l} V_{3} & = & \ds{(x^{2}-3y^{2})\frac{\partial}{\partial x} + (t+4xy)\frac
{\partial}{\partial y}}\\ \\
& &  \ds{+ (2xt-2x^{2}y-2y^{3})\frac{\partial}{\partial t} - 2x
u \frac{\partial}{\partial u}},
\ea\end{equation}
\[Z =x\frac{\partial}{\partial x}+y\frac{\partial}{\partial y}+2t
\frac{\partial}{\partial t},\]
\bb\label{uw} U= u\frac{\partial}{\partial
u},\;\;W_{\beta}= \beta(x,y,t)\frac{\partial}{\partial u}, \text{ where }
\Delta_{H^{1}} \beta=0.
\end{equation}
{\footnotesize
\begin{table}[h]
\begin{center}%
\begin{tabular}
[c]{|c|c|c|c|c|c|c|c|c|c|c|}\hline
& T & R & $\tilde{X}$ & $\tilde{Y}$ & U & $W_{\beta}$ & $V_{1}$ & $V_{2}$ &
$V_{3}$ & $Z$\\\hline
T & $0$ & $0$ & $0$ & $0$ & $0$ & $W_{T\beta}$ & $V$ & $\tilde{X}$ &
$\tilde{Y}$ & $2$T\\\hline
R & $0$ & $0$ & $\tilde{Y}$ & -$\tilde{X}$ & $0$ & $W_{R\beta}$ & $0$ &
$V_{3}$ & -$V_{2}$ & $0$\\\hline
$\tilde{X}$ & $0$ & -$\tilde{Y}$ & $0$ & $4$T & $0$ & $W_{\tilde{X}\beta}$ &
$V_{2}$ & -$6$R & $2V$ & $\tilde{X}$\\\hline
$\tilde{Y} $ & $0$ & $\tilde{X}$ & -$4T$ & $0$ & $0$ & $W_{\tilde{Y}\beta}$ &
$V_{3}$ & -$2V$ & -$6$R & $\tilde{Y}$\\\hline
U & $0$ & $0$ & $0$ & $0$ & $0$ & $0$ & $0$ & $0$ & $0$ & $0$\\\hline
$W_{\beta}$ & -$W_{T\beta}$ & -$W_{R\beta}$ & -$W_{\tilde{X}\beta}$ &
-$W_{\tilde{Y}\beta}$ & $0$ & 0 & $W_{V_{1}\beta}$ & $W_{V_{2}\beta}$ &
$W_{V_{3}\beta}$ & $W_{Z\beta}$\\\hline
$V_{1}$ & -$V$ & $0$ & -$V_{2}$ & -$V_{3}$ & $0$ & -$W_{V_{1}\beta}$ & $0$ &
$0$ & $0$ & -$2V_{1}$\\\hline
$V_{2}$ & -$\tilde{X}$ & -$V_{3}$ & $6R$ & $2V$ & $0$ & -$W_{V_{2}\beta}$ &
$0$ & $0$ & $4 V_{1}$ & -$V_{2}$\\\hline
$V_{3}$ & -$\tilde{Y}$ & $V_{2}$ & -$2V$ & $6R$ & $0$ & -$W_{V_{3}\beta}$ &
$0 $ & -$4$$V_{1}$ & $0$ & -$V_{3}$\\\hline
$Z$ & -$2$T & $0$ & -$\tilde{X}$ & -$\tilde{Y}$ & $0$ & -$W_{Z\beta}$ &
$2V_{1}$ & $V_{2}$ & $V_{3}$ & $0$\\\hline
\end{tabular}
\end{center}
\caption{{\small {Lie brackets of equation $(\ref{klgen})$ with $f(u)=0$}. Here, $V:=Z-U$.}}%
\end{table}}

\item If $f(u)=u $, there are two additional symmetries, respectively, $U$ and $W_{\beta}$ as in Eq. (\ref{uw}), where $\Delta_{H^{1}} \beta
+\beta=0$.
\begin{table}[h]
\begin{center}%
\begin{tabular}
[c]{|c|c|c|c|c|c|c|}\hline
& T & R & $\tilde{X}$ & $\tilde{Y}$ & U & $W_{\beta}$\\\hline
T & $0$ & $0$ & $0$ & $0$ & $0$ & $W_{T\beta}$\\\hline
R & $0$ & $0$ & $\tilde{Y}$ & -$\tilde{X}$ & $0$ & $W_{R\beta}$\\\hline
$\tilde{X}$ & $0$ & -$\tilde{Y}$ & $0$ & $4$T & $0$ & $W_{\tilde{X}\beta}%
$\\\hline
$\tilde{Y}$ & $0$ & $\tilde{X}$ & $4$T & $0$ & $0$ & $W_{\tilde{Y}\beta}%
$\\\hline
U & $0$ & $0$ & $0$ & $0$ & $0$ & $0$\\\hline
$W_{\beta}$ & -$W_{T\beta}$ & -$W_{R\beta}$ & -$W_{\tilde{X}\beta}$ &
-$W_{\tilde{Y}\beta}$ & $0$ & $0$\\\hline
\end{tabular}
\end{center}
\caption{{\small {Lie brackets of equation $(\ref{klgen})$ with $f(u)=u$}.}}%
\end{table}

\item If $f(u)=u^{p} $, $p\neq0, p\neq1, p\neq3$, we have the generator of
dilations
\begin{equation}
\label{D}D_{p}= x\frac{\partial}{\partial x}+y\frac{\partial}{\partial y}+2t
\frac{\partial}{\partial t}+\frac{2}{1-p}u\frac{\partial}{\partial u}.
\end{equation}
\begin{table}[h]
\begin{center}%
\begin{tabular}
[c]{|c|c|c|c|c|c|}\hline
& T & R & $\tilde{X}$ & $\tilde{Y}$ & $D_{p}$\\\hline
T & $0$ & $0$ & $0$ & $0$ & $2$T\\\hline
R & $0$ & $0$ & $\tilde{Y}$ & -$\tilde{X}$ & $0$\\\hline
$\tilde{X}$ & $0$ & -$\tilde{Y}$ & $0$ & $4$T & $\tilde{X}$\\\hline
$\tilde{Y} $ & $0$ & $\tilde{X}$ & - $4T$ & $0$ & $\tilde{Y}$\\\hline
$D_{p}$ & -$2$T & $0$ & -$\tilde{X}$ & -$\tilde{Y}$ & $0$\\\hline
\end{tabular}
\end{center}
\caption{{\small {Lie brackets of equation $(\ref{klgen})$ with $f(u)=u^{p}$,
$p\neq0, p\neq1, p\neq3$.}}}%
\end{table}

\item If $f(u)=e^{u}$ the additional symmetry is
\begin{equation}
\label{Z}E=x\frac{\partial}{\partial x}+y\frac{\partial}{\partial y}+2t
\frac{\partial}{\partial t}-2\frac{\partial}{\partial u}.
\end{equation}
\begin{table}[h]
\begin{center}%
\begin{tabular}
[c]{|c|c|c|c|c|c|}\hline
& T & R & $\tilde{X}$ & $\tilde{Y}$ & $E$\\\hline
T & $0$ & $0$ & $0$ & $0$ & $2$T\\\hline
R & $0$ & $0$ & $\tilde{Y}$ & -$\tilde{X}$ & $0$\\\hline
$\tilde{X}$ & $0$ & -$\tilde{Y}$ & $0$ & $4$T & $\tilde{X}$\\\hline
$\tilde{Y} $ & $0$ & $\tilde{X}$ & - $4$T & $0$ & $\tilde{Y}$\\\hline
$E$ & -$2$T & $0$ & -$\tilde{X}$ & -$\tilde{Y}$ & $0$\\\hline
\end{tabular}
\end{center}
\caption{{\small {Lie brackets of equation $(\ref{klgen})$ with $f(u)=e^{u}$%
}.}}%
\end{table}

\item In the critical case, $f(u)=u^{ 3}$, there are four additional
generators, namely $V_{1},V_{2},V_{3}\text{ and } D_{3}$, given in
$(\ref{V1}), (\ref{V2}),(\ref{V3})\text{ and }(\ref{D})$ respectively. Their
Lie algebra is presented in \cite{cl}.
\end{itemize}

In \cite{ds} is showed that in the critical case, $f(u)=u^{3}$, all Lie
point symmetries are Noether symmetries and then, by the Noether Theorem (see
\cite{bl}, pag. 275), in \cite{cl} is established the respectives conservation
laws for the symmetries $T,R,\tilde{X},\tilde{Y},V_{1},V_{2},V_{3}$ and
$D_{3}$.

In this work, we show which Lie point symmetries of the other functions $f(u)$
are Noether symmetries and then, we establish their respectives conservation
laws, concluding the work started in \cite{ds} and \cite{cl}.

Let $\mathbb{R}\ni u\mapsto F(u)\in\mathbb{R}$ be a differentiable function
and
\begin{equation}
\label{F}f(u):=F^{\prime}(u).
\end{equation}
Our main results can be formulated as follows:

\begin{theorem}
\label{noethergeral} The group $G_{f}$ is a Noether symmetry group for any
function $f(u)$ in $(\ref{klgen})$.
\end{theorem}

\begin{theorem}
\label{noetherexp} The Noether symmetry group of $(\ref{klgen})$, with
$f(u)=e^{u}$, is the group $G_{f}$.
\end{theorem}

\begin{theorem}
\label{noetl} $G_{f}\cup\{W_{\beta}\}$ is the Noether symmetry group of
equation $(\ref{klgen})$, with $f(u)=u$ and $\beta$ satisfies $\Delta_{H^{1}%
}\beta+\beta=0$.
\end{theorem}

\begin{theorem}
\label{noethom} The Noether symmetry group of equation $(\ref{klgen})$ with
$f(u)=0$ is generated by the group $G_{f}$ and by symmetries $W_{\beta
},\;V_{1},\;V_{2}$ e $V_{3}$, where $\beta$ satisfies $\Delta_{H^{1}}\beta=0$.
\end{theorem}

As a consequence of theorems \ref{noethergeral} - \ref{noethom}, we have the
following conservation laws.

\begin{theorem}
\label{leigeral} The conservations laws for the Noether symmetries of equation
$(\ref{klgen})$ for any $f(u)$ are:

\begin{enumerate}
\item For the symmetry $T$, the conservation law is $Div(\tau)=0$, where
$\tau=(\tau_{1},\tau_{2},\tau_{3})$ and
\[%
\begin{array}
[c]{l}%
\ds{\tau_{1}=-2}yu{_{t}^{2}-}u{_{x}}u{_{t}},\\

\ds{\tau_{2}=2}xu{_{t}^{2}-}u_{y}u{_{t}},\\

\ds{\tau_{3}=\frac{1}{2}}u{_{x}^{2}+\frac{1}{2}}u{_{y}^{2}-2(}x{^{2}+}y{^{2}%
)}u{_{t}^{2}-F(}u{)}.
\end{array}
\]

\item For the symmetry $R$, the conservation law is $Div(\sigma)=0$, where
$\sigma=(\sigma_{1},\sigma_{2},\sigma_{3})$ and
\[%
\begin{array}
[c]{l}%
\ds{\sigma_{1}=-\frac{1}{2}yu_{x}^{2}+\frac{1}{2}yu_{y}^{2}+2y(x^{2}+y^{2}%
)u_{t}^{2}+xu_{x}u_{y}-yF(u)},\\

\ds{\sigma_{2}=-\frac{1}{2}xu_{x}^{2}-\frac{1}{2}xu_{y}^{2}-2x(x^{2}+y^{2}%
)u_{t}^{2}-yu_{x}u_{y}+xF(u)},\\

\ds{\sigma_{3}=-2y^{2}u_{x}^{2}-2x^{2}u_{y}^{2}+4xyu_{x}u_{y}-4y(x^{2}
+y^{2})u_{x}u_{t}+4x(x^{2}+y^{2})u_{y}u_{t}}.
\end{array}
\]

\item For the symmetry $\tilde{X}$, the conservation law is $Div(\chi)=0$,
where $\chi=(\chi_{1},\chi_{2},\chi_{3})$ and
\[%
\begin{array}
[c]{l}%
\ds{\chi_{1}=}-\frac{1}{2}u_{x}^{2}+\frac{1}{2}u_{y}^{2}+2(x^{2}+3y^{2})u_{t}
^{2}+2yu_{x}u_{t}-2xu_{y}u_{t}-F(u),\\

\ds{\chi_{2}=}-4xyu_{t}^{2}-u_{x}u_{y}+2xu_{x}u_{t}+2yu_{y}u_{t},\\

\ds{\chi_{3}=-3yu_{x}^{2}-yu_{y}^{2}+4y(x^{2}+y^{2})u_{t}^{2}+2xu_{x}
u_{y}-4(x^{2}+y^{2})u_{x}u_{t}+2yF(u)}.
\end{array}
\]

\item For the symmetry $\tilde{Y}$, the conservation law is $Div(\upsilon)=0$,
where $\upsilon=(\upsilon_{1},\upsilon_{2},\upsilon_{3})$ and
\[%
\begin{array}
[c]{l}%
\ds{\upsilon_{1}=-4xyu_{t}^{2}-u_{x}u_{y}-2xu_{x}u_{t}-2yu_{y}u_{t}},\\ \\

\ds{\upsilon_{2}=\frac{1}{2}u_{x}^{2}-\frac{1}{2}u_{y}^{2}+2(3x^{2}+y^{2}
)u_{t}^{2}+2yu_{x}u_{t}-2xu_{y}u_{t}-F(u)},\\ \\

\ds{\upsilon_{3}=xu_{x}^{2}+3xu_{y}^{2}-4x(x^{2}+y^{2})u_{t}^{2}-2yu_{x}
u_{y}-4(x^{2}+y^{2})u_{y}u_{t}-2xF(u)}.
\end{array}
\]

\end{enumerate}
\end{theorem}

\begin{theorem}
\label{leihomogen} If $f(u)=0$ in $(\ref{klgen})$, the conservation laws for
the Noether symmetries are as follows.

\begin{enumerate}
\item For the symmetries $T,\;R,\;\tilde{X}$ and $\tilde{Y}$, the conservation
laws are the same as in the Theorem $\ref{leigeral}$, with $f(u)=0$, in
$(\ref{F})$.

\item For the symmetry $V_{1}$, the conservation law is $Div(A)=0$, where
$A=(A_{1},A_{2},A_{3})$ and%
\begin{align*}
A_{1}  & \ds{=-\frac{1}{2}(tx-x^{2}y-y^{3})u_{x}^{2}+\frac{1}{2}(tx-x^{2}
y-y^{3})u_{y}^{2}+2t(x^{3}+xy^{2}-ty)u_{t}^{2}}\\ \\
& \ds{-(x^{3}+xy^{2}+ty)u_{x}u_{y}-[t^{2}-(x^{2}+y^{2})^{2}]u_{x}u_{t}
-2t(x^{2}+y^{2})u_{y}u_{t}}\\ \\
& \ds{-tuu_{x}-2tyuu_{t}+yu^{2}},
\end{align*}%
\begin{align*}
A_{2}  & \ds{=\frac{1}{2}(x^{3}+ty+xy^{2})u_{x}^{2}-\frac{1}{2}(x^{3}
+ty+xy^{2})u_{y}^{2}+2t(x^{2}y+y^{3}+tx)u_{t}^{2}}\\ \\
& \ds{-(tx-x^{2}y-y^{3})u_{x}u_{y}+2t(x^{2}+y^{2})u_{x}u_{t}-[t^{2}-(x^{2}
+y^{2})^{2}]u_{y}u_{t}}\\ \\
&\ds{ -tuu_{y}+2txuu_{t}-xu^{2}},
\end{align*}%
\begin{align*}
A_{3}  & \ds{=\frac{1}{2}(t^{2}-x^{4}-4txy+2x^{2}y^{2}+3y^{4})u_{x}^{2}+\frac
{1}{2}(t^{2}+3x^{4}+4txy+2x^{2}y^{2}-y^{4})u_{y}^{2}}\\ \\
& \ds{-2(x^{2}+y^{2})[t^{2}-(x^{2}+y^{2})^{2}]u_{t}^{2}+2[t(x^{2}-y^{2}
)-2xy(x^{2}+y^{2})]u_{x}u_{y}}\\ \\
& \ds{-4(x^{2}+y^{2})(tx-x^{2}y-y^{3})u_{x}u_{t}-4(x^{2}+y^{2})(x^{3}
+ty+xy^{2})u_{y}u_{t}}\\ \\
& \ds{-2tyuu_{x}+2txuu_{y}-4t(x^{2}+y^{2})uu_{t}+2(x^{2}+y^{2})u^{2}}.
\end{align*}

\item For the symmetry $V_{2}$, the conservation law is $Div(B)=0$, where
$B=(B_{1},B_{2},B_{3})$ and%
\begin{align*}
B_{1}  &\ds{ =-\frac{1}{2}(t-4xy)u_{x}^{2}+\frac{1}{2}(t-4xy)u_{y}^{2}
+[2t(x^{2}+3y^{2})-4xy(x^{2}+y^{2})]u_{t}^{2}}\\ \\
& \ds{-(3x^{2}-y^{2})u_{x}u_{y}+2(x^{3}+ty+xy^{2})u_{x}u_{t}-2(tx-x^{2}
y-y^{3})u_{y}u_{t}}\\ \\
&\ds{ +2yuu_{x}+4y^{2}uu_{t}},
\end{align*}%
\begin{align*}
B_{2}  & \ds{=\frac{1}{2}(3x^{2}-y^{2})u_{x}^{2}-\frac{1}{2}(3x^{2}-y^{2}
)u_{y}^{2}+2(x^{4}-2txy-y^{4})u_{t}^{2}-(t-4xy)u_{x}u_{y}}\\ \\
& \ds{+2(tx-x^{2}y-y^{3})u_{x}u_{t}+2(x^{3}+ty+xy^{2})u_{y}u_{t}+2yuu_{y}
-4xyuu_{t}-u^{2}},
\end{align*}%
\begin{align*}
B_{3}  & \ds{=(7xy^{2}-x^{3}-3ty)u_{x}^{2}+(5x^{3}-3xy^{2}-ty)u_{y}^{2}
+4(x^{2}+y^{2})(x^{3}+ty+xy^{2})u_{t}^{2}}\\ \\
& \ds{+2(tx-7x^{2}y+y^{3})u_{x}u_{y}-4(t-4xy)(x^{2}+y^{2})u_{x}u_{t}
-4(3x^{4}+2x^{2}y^{2}-y^{4})u_{y}u_{t}}\\ \\
& \ds{+2xu^{2}+4y^{2}uu_{x}-4xyuu_{y}+8y(x^{2}+y^{2})uu_{t}}.
\end{align*}

\item For the symmetry $V_{3}$, the conservation law is $Div(C)=0$, where
$C=(C_{1},C_{2},C_{3})$ and%
\begin{align*}
C_{1}  & \ds{=-\frac{1}{2}(x^{2}-3y^{2})u_{x}^{2}+\frac{1}{2}(x^{2}-3y^{2}
)u_{y}^{2}+(2x^{4}-4txy-2y^{4})u_{t}^{2}}\\ \\
& \ds{-(t+4xy)u_{x}u_{y}+(2tx-2x^{2}y+2y^{3})u_{x}u_{t}-(2x^{3}+2ty+2xy^{2}
)u_{y}u_{t}}\\ \\
& \ds{-4xyuu_{t}-2xuu_{x}+u^{2}},
\end{align*}%
\begin{align*}
C_{2}  & \ds{=\frac{1}{2}(t+4xy)u_{x}^{2}-\frac{1}{2}(t+4xy)u_{y}^{2}
+(6tx^{2}+4x^{3}y+2ty^{2}+4xy^{3})u_{t}^{2}}\\ \\
& \ds{-(x^{2}-3y^{2})u_{x}u_{y}+2(x^{3}+ty+xy^{2})u_{x}u_{t}-2(tx-x^{2}
y-y^{3})u_{y}u_{t}}\\ \\
& \ds{2xu_{y}u+4x^{2}u_{t}u},
\end{align*}%
\begin{align*}
C_{3}  & \ds{=(tx-3x^{2}y+5y^{3})u_{x}^{2}+(3tx+7x^{2}y-y^{3})u_{y}^{2}}\\ \\
& \ds{(-4tx^{3}+4x^{4}y-4txy^{2}+8x^{2}y^{3}+y^{5})u_{t}^{2}+2(x^{3}-ty-7xy^{2}
)u_{x}u_{y}}\\ \\
& \ds{-2(2x^{4}-4x^{2}y^{2}-6y^{4})u_{x}u_{t}-4(x^{2}+y^{2})(t+4xy)u_{y}u_{t}}\\ \\
& \ds{-8x^{3}uu_{t}-8xy^{2}uu_{t}-4x^{2}uu_{y}-8xyuu_{x}+2yu^{2}}.
\end{align*}

\item For the symmetry $W_{\beta}$, the conservation law is $Div(W)=0$, where
$W=(W_{1},W_{2},W_{3})$ and
\begin{equation}\label{leiw}\begin{array}{l c l}
W_{1} & = & \ds{\beta(u_{x}+2yu_{t})-u(\beta_{x}+2y\beta_{t})},\\ \\

W_{2} & = & \ds{\beta(u_{y}-2xu_{t})-u(\beta_{y}-2x\beta_{t})},\\ \\

W_{3} & = & \ds{\beta\lbrack-2xu_{y}+2yu_{x}+4(x^{2}+y^{2})u_{t}}\\ \\

& & \ds{+2u[x\beta
_{y}-y\beta_{x}-2(x^{2}+y^{2})\beta_{t}]}.
\end{array}\end{equation}
\end{enumerate}
\end{theorem}

\begin{theorem}
\label{leilin} If $f(u)=u$ in $(\ref{klgen})$, the conservation laws for the
Noether symmetries are as follows.

\begin{enumerate}
\item For the symmetries $T,\;R,\;\tilde{X}$ and $\tilde{Y}$, the conservation
laws are the same as in the Theorem $\ref{leigeral}$, with $f(u)=u$, in
$\emph{(\ref{F})}$.

\item For the symmetry $W_{\beta}$, the conservation law is $Div(W)=0$, where
$W$ is given in $\ref{leiw}$ .
\end{enumerate}
\end{theorem}

The remaining of the paper is organized as follows. In section 2 we briefly
present some of the main aspects of Lie point symmetries, Noether symmetries
and conservation laws. In section 3 we prove theorems \ref{noethergeral},
\ref{noetherexp} and \ref{noetl}. Theorem \ref{noethom} is proved in section
4. Their respective conservation laws are discussed in section 5.

\section{Lie point symmetries, Noether symmetries and conservation laws}

Let $x\in M\subseteq\mathbb{R}^{n}$, $u:M\rightarrow\mathbb{R}$ a smooth
function and $k\in\mathbb{N}$. $\partial^{k} u$ denotes the jet bundle
correspondig to all $k$th partial derivatives of $u$ with respect to $x$. A
\textit{Lie point symmetry} of a partial differential equation (PDE) of order
$k$, $F(x,u,\partial u,\cdots,\partial^{k} u)=0$, is a vector field
\[
S=\xi^{i}(x,u)\frac{\partial}{\partial x^{i}}+\eta(x,u)\frac{\partial
}{\partial u}%
\]
on $M\times\mathbb{R}$ such that $S^{(k)}F=0$ when $F=0$ and
\[
S^{(k)}=S+\eta^{(1)}_{i}(x,u,\partial u)\frac{\partial}{\partial u_{i}}%
+\cdots+\eta^{(k)}_{i_{1}\cdots i_{k}}(x,u,\partial u,\cdots,\partial^{k}
u)\frac{\partial}{\partial u_{i_{1}\cdots i_{k}}}%
\]
is the extended symmetry on the jet space $(x,u,\partial u,\cdots,
\partial^{k} u)$.

The functions $\eta^{(j)}(x,u,\partial u,\cdots,\partial^{j} u)$, $1\leq j\leq
k$, are given by
\[
\begin{array}
[c]{lcl}%
\eta^{(1)}_{ i} & = & D_{i}\eta-(D_{i}\xi^{j})u_{j},\\ \\
&  & \\ \\
\eta^{(j)}_{i_{1}\cdots i_{j}} & = & D_{i_{j}}\eta^{(j-1)}_{i_{1}\cdots
i_{j-1}}-(D_{i_{j}}\xi^{l})u_{i_{1}\cdots i_{j-1}l},\;2\leq j\leq k.
\end{array}
\]
We are using the Einstein sum convention.

If the PDE $F=0$ can be obtained by a Lagrangian $\mathcal{L}=\mathcal{L}%
(x,u,\partial u,\cdots, \partial^{l} u)$ and if there exists some symmetry $S$
of $F$ and a vector $\varphi=(\varphi_{1},\cdots,\varphi_{n})$ such that
\begin{equation}
\label{defno}S^{(l)}\mathcal{L}+\mathcal{L}D_{i}\xi^{i}=D_{i}\varphi^{i},
\end{equation}
where
\[
D_{i}=\frac{\partial}{\partial x^{i}}+u_{i}\frac{\partial}{\partial u}%
+u_{ij}\frac{\partial}{\partial u_{j}}+\cdots+u_{ii_{1}\cdots i_{m}}%
\frac{\partial}{\partial u_{i_{1}\cdots i_{m}}}+\cdots
\]
is the total derivative operator of $u$,
\[
u_{i}:=\frac{\partial u}{\partial x^{i}},\; u_{ij}:=\frac{\partial^{2}%
u}{\partial x^{i}\partial x^{j}},\cdots, u_{ii_{1}\cdots i_{m}}:=\frac
{\partial u}{\partial x_{i}\partial x_{i_{1}}\cdots\partial x_{i_{m}}},
\]
the symmetry $S$ is said to be a \textit{Noether symmetry}. Then, the
Noether's Theorem asserts that the following conservation law holds
\begin{equation}
\label{cons}D_{i}(\xi^{i}\mathcal{L}+W^{i}[u,\eta-\xi^{j}u_{j}]-\varphi
^{i})=0.
\end{equation}
Above we have used the same notations and conventions as in \cite{bl}. (For
the definition of $W^{i}$ see \cite{bl}, pp. 254-255.)

\section{Proofs of theorems 1, 2 and 3}

\begin{lemma}
\label{tinexiste} Let $u=u(x,y,t)$ be a smooth function. If a vector field
$V=(A,B,C)$ is a vector function of $x$, $y$, $t$, $u$, $u_{x}$, $u_{y}$,
$u_{t}$, its divergence necessarily depends on the second order derivatives of
u with respect to $x,\;y$ and $t$.
\end{lemma}

\begin{proof}
Taking the divergence of vector field $V$, we obtain
\begin{align*}
Div(V)=  & A_{x}+B_{y}+C_{t}+u_{x}A_{u}+u_{xx}A_{u_{x}}+u_{xy}A_{u_{y}}%
+u_{xt}A_{u_{t}}\\ \\
& +u_{y}B_{u}+u_{xy}B_{u_{x}}+u_{yy}B_{u_{y}}+u_{yt}B_{u_{t}}\\ \\
& +u_{t}C_{u}+u_{xt}C_{u_{x}}+u_{yt}C_{u_{y}}+u_{tt}C_{u_{t}}.
\end{align*}
\end{proof}
\begin{corollary}
\label{cinexiste} If the divergence of a vector field does not depend on the
second order derivatives, then it does not depend on $u_{x},\;u_{y}$ and
$u_{t}$.
\end{corollary}
\begin{lemma}
\label{lemmma} The symmetry
\[
E=x\frac{\partial}{\partial x}+y\frac{\partial}{\partial y}+2t\frac{\partial
}{\partial t}-2\frac{\partial}{\partial u}%
\]
is not a Noether symmetry.
\end{lemma}
\begin{proof}
In this case, $(\xi,\phi,\tau,\eta)=(x,y,2t,-2)$. Then, $D_{x}\xi+D_{y}%
\phi+D_{t}\tau=4$ and
\[
(\eta^{(1)}_{x},\eta^{(1)}_{y},\eta^{(1)}_{t})=(-u_{x},-u_{y},-2u_{t}),
\]
which yields the following first order extension:
\[
E^{(1)}=E-u_{x}\f{\p}{\p u_{x}}-u_{y}\f{\p}{\p
u_{y}}-2u_{t}\f{\p}{\p u_{t}}.
\]
Therefore, 
\begin{equation}
\label{extE}
\ba{l c l}
\ds{E^{(1)}\mathcal{L}+(D_{x}\xi+D_{y}\phi+D_{t}\tau)\cal{L}}& = &
\ds{u_{x}^{2}+u_{y}^{2}+4(x^{2}+y^{2})u_{t}^{2}}
\\ \\
&
&\ds{+4yu_{x}u_{t}-4xu_{y}u_{t}-2e^{u}},
\ea
\end{equation}
where
\[
\mathcal{L}:=\frac{1}{2}u_{x}^{2}+\frac{1}{2}u_{y}^{2}+2(x^{2}+y^{2})u_{t}%
^{2}+2yu_{x}u_{t}-2xu_{y}u_{t}-e^{u}.
\]
From Lemma \ref{tinexiste} and equation (\ref{extE}), we conclude that there
are not a potential $\phi$ which satisfies
\[
E^{(1)}\mathcal{L}+(D_{x}\xi+D_{y}\phi+D_{t}\tau)\mathcal{L}=Div(\phi).
\]
Thus, $E$ cannot be a Noether symmetry.
\end{proof}
\begin{lemma}
\label{lemau} The symmetry $U$ is not a Noether symmetry.
\end{lemma}
\begin{proof}
First one, note that $\eta=u,\;\xi=\phi=\tau=0$. Then,
\begin{equation}
U^{(1)}=u\f{\p}{\p
u}+u_{x}\f{\p}{\p u_{x}}+u_{y}\f{\p}{\p u_{t}}+u_{t}\f{\p}{\p
u_{t}}\label{extu}%
\end{equation}
Aplying the operator obtained in (\ref{extu}) to the Lagrangian
\begin{equation}
\mathcal{L}_{k}:=\frac{1}{2}u_{x}^{2}+\frac{1}{2}u_{y}^{2}+2(x^{2}+y^{2}%
)u_{t}^{2}+2yu_{x}u_{t}-2xu_{y}u_{t}-\f{k}{2}u^{2},\label{lagl}%
\end{equation}
where $\ k=0$ if $f(u)=0$ or $k=1$ if $f(u)=u$, we find
\[
U^{(1)}\mathcal{L}_{k}=u_{x}^{2}+u_{y}^{2}+4(x^{2}+y^{2})u_{t}^{2}%
+4yu_{x}u_{t}-4xu_{y}u_{t}-ku^{2}=2\mathcal{L}_{k}.
\]
From Lemma \ref{tinexiste} and Corollary \ref{cinexiste}, we conclude that
there is not a vector field such that equation (\ref{defno}) is true with
$S=U$.
\end{proof}
\begin{lemma}
\label{lemaw} The symmetry $W_{\beta}$ is a Noether symmetry.
\end{lemma}
\begin{proof}
The first order extension $W^{(1)}$ of $W$ is
\begin{equation}
W^{(1)}=\be\f{\p}{\p u}+\be_{x}\f{\p}{\p
u_{x}}+\be_{y}\f{\p}{\p u_{y}}+\be_{t}\f{\p}{\p u_{t}}.\label{w}%
\end{equation}
From (\ref{w}) and (\ref{lagl}), we have
\begin{align*}
W^{(1)}\mathcal{L}_{k}  & \ds{=-\be ku+(u_{x}+2yu_{t})\be_{x}+(u_{y}
-2xu_{t})\be_{y}+(4(x^{2}+y^{2})u_{t}+2yu_{x}-2xu_{y})\be_{t}}\\ \\
\\ \\
& \ds{=Div((\be_{x}+2y\be_{t})u,(\be_{y}-2x\be_{t})u,(2y\be_{x}-2x\be_{y}
+4(x^{2}+y^{2})\be_{t})u)}.
\end{align*}
\end{proof}
\begin{lemma}
\label{lemaZ} The symmetry
\[
Z =x\frac{\partial}{\partial x}+y\frac{\partial}{\partial y}+2t \frac
{\partial}{\partial t}%
\]
is not a Noether symmetry.
\end{lemma}
\begin{proof}
Since $D_{x}\xi+D_{y}\phi+D_{t}\tau=4$,
\begin{equation}
\mathcal{L}=\frac{1}{2}u_{x}^{2}+\frac{1}{2}u_{y}^{2}+2(x^{2}+y^{2})u_{t}%
^{2}+2yu_{x}u_{t}-2xu_{y}u_{t}\label{extZ}%
\end{equation}
\ and%
\begin{equation}
Z^{(1)}=Z+u_{x}\f{\p}{\p x}+u_{y}\f{\p}{\p
y}+2u_{t}\f{\p}{\p t}\label{lagh}%
\end{equation}
is a consequence of Eqs. (\ref{extZ}) and (\ref{lagh}), that
\begin{equation}\label{ver}
Z^{(1)}\mathcal{L}+\mathcal{L}(D_{x}\xi+D_{y}\phi+D_{t}\tau)=2\mathcal{L}.
\end{equation}
By Lemma \ref{tinexiste}, there does not exist a vector field such that the
right hand of (\ref{ver}) be its divergence.
\end{proof}
\textbf{Proof of Theorem \ref{noethergeral}}: We will use four steps to prove
this theorem. First, we obtain the first order extension of symmetries
$T,\;R,\;\tilde{X},\;\tilde{Y}$. Next, we proof the theorem for each one of them.
\begin{enumerate}
\item Extensions:
\begin{enumerate}
\item Symmetry $T$
The coefficients of $T$ are $\xi=\phi=\eta=0$ and $\phi=1.$ Then
\[
T^{(1)}=T.
\]
\item Symmetry $R$
The coefficients of symmetry $R$ are $(\xi,\phi,\tau,\eta)=(y,-x,0,0) $. Then,
we conclude that
\[
R^{(1)}=R+u_{y}\frac{\partial}{\partial u_{x}}- u_{x}\frac{\partial}{\partial
u_{y}}.
\]
\item Symmetry $\tilde{X}$
In this case, $(\xi,\phi,\tau,\eta)=(1,0,-2y,0,)$. Then
\[%
\begin{array}
[c]{lcl}%
\eta^{(1)}_{x}=0, & \eta_{y}^{(1)}=2u_{t}, & \eta^{(1)}_{t}=0
\end{array}
\]
and
\[
{\tilde{X}}^{(1)}=\tilde{X} + 2u_{t}\frac{\partial}{\partial u_{y}}.
\]
\item Symmetry $\tilde{Y}$
This case is analogous to case $c$ and we present only its extension
\[
{\tilde{Y}}^{(1)}=\tilde{Y} - 2u_{t}\frac{\partial}{\partial u_{x}}.
\]
\end{enumerate}
\begin{corollary}
\label{divk} The divergence of any symmetry $S\in G_{f}$ is zero.
\end{corollary}
\item
\begin{enumerate}
\item {Proof of theorem for the symmetry $T$}.
Since $Div(T)=0= T^{(1)}\mathcal{L}$ it is immediate that
\[
T^{(1)}\mathcal{L}+\mathcal{L}Div(T)=0.
\]
\item {Proof of theorem for the symmetry $R$}.
We have
\[
R^{(1)}\mathcal{L}=0.
\]
Then, from Corollary \ref{divk},
\[
R^{(1)}\mathcal{L}+\mathcal{L}Div(R)=0.
\]
\item {Proof of theorem for the symmetries $\tilde{X}$ and $\tilde{Y}$}.
It is immediate that
\[
\tilde{X}^{(1)}\mathcal{L}=0.
\]
Again, by Corollary \ref{divk}, we obtain
\[
\tilde{X}^{(1)}\mathcal{L}+\mathcal{L}Div(\tilde{X})=0.
\]
In the same way, we conclude that
\[
\tilde{Y}^{(1)}\mathcal{L}+\mathcal{L}Div(\tilde{Y})=0.
\]
\end{enumerate}
\end{enumerate}
\textbf{Proof of Theorem \ref{noetherexp}}: It is a consequence of Lemma
\ref{lemmma} and Theorem \ref{noethergeral}.\newline
\textbf{Proof of Theorem \ref{noetl}}: From Lemma \ref{lemau}, $U$ is not a
Noether symmetry. Then, by Theorem \ref{noethergeral} and Lemma \ref{lemaw},
$G_{f}\cup\{W_{\beta}\}$ is a Noether symmetry group.\newline
\textbf{Proof of Theorem \ref{noethom}}: By lemmas \ref{lemau} and
\ref{lemaZ}, the symmetries $Z$ and $U$ are not Noether symmetries. The proof
that the symmetries $V_{1},\;V_{2}$ and $V_{3}$ are Noether symmetries is
obtained in same way that Bozhkov and Freire showed that $V_{1},\;V_{2}$ and
$V_{3}$ are Noether symmetries of (\ref{klgen}) when $f(u)=u^{3}$, and can be
found in \cite{ds}. Then, by Theorem \ref{noethergeral} and Lemma \ref{lemaw},
we conclude the proof.

\section{Conservation Laws}

The proof is by a straightforward calculation, which we shall not present
here. However, a computer assisted proof can be obtained by means of the
software \textit{Mathematica}. It calculates the components of the
conservation laws, which appear in the equation (\ref{cons}). The Mathematica
notebook used for this purpose can be obtained form the author upon request.

\begin{center}
\textbf{Acknowledgements}
\end{center}

We are grateful to Yuri Bozhkov for his suggestions and firm encouragement to
write this paper. We also show gratitude to Antonio Carlos Gilli Martins,
Ricardo Antonio Mosna and Waldyr Alves Rodrigues Jr. for their suggestions and helpful.

We also thank Lab. EPIFISMA (Proj. FAPESP) for having given us the opportunity
to use excellent computer facilities and we are grateful to UNICAMP for
financial support.

\end{document}